\documentclass[12pt,amsfonts, epsfig]{amsart}
\usepackage{amsmath, amscd, amssymb}

\usepackage[frame,cmtip,arrow,matrix,line,graph,curve]{xy}
\usepackage{epsfig}
\usepackage{graphpap, color}
\usepackage[mathscr]{eucal}
\usepackage{mathrsfs}

\numberwithin{equation}{section}

\newcommand{\Gr}{\operatorname{Gr}}

\newcommand{\PP}{\mathbb{P}}

\newcommand{\bR}{\mathbf{R}}

\newcommand{\bk}{\mathbf{k}}

\newcommand{\kk}{\bk}

\newcommand{\cal}{\mathcal}

\def\cA{{\cal A}}

\def\cB{{\cal B}}
\def\cC{{\cal C}}
\def\cD{{\cal D}}

\def\cE{{\cal E}}

\def\cL{{\cal L}}

\def\cO{{\cal O}}
\def\cP{{\cal P}}
\def\cQ{{\cal Q}}

\def\cV{{\cal V}}

\def\cX{{\cal X}}
\def\cY{{\cal Y}}

\def\Eb{E^\bullet}

\def\sE{{\mathscr E}}
\def\sEb{\sE^\bullet}

\def\sH{{\mathscr H}}
\def\sI{{\mathscr I}}

\def\sL{{\mathscr L}}

\def\sO{{\mathscr O}}

\def\fM{\mathfrak{M}}
\def\2M{M}

\def\fP{\mathfrak{P}}

\def\fX{\mathfrak{X}}

\def\ff{\mathfrak{f}}

\def\dual{^{\vee}}

\def\sta{^\ast}

\def\virt{^{\mathrm{vir}}}

\def\sta{^{\ast}}

\def\sta{^*}

\def\lra{\longrightarrow}

\newcommand{\Ga}{\Gamma}

\newcommand{\ee}{{\bf e}}

\def\begeq{\begin{equation}}
\def\endeq{\end{equation}}
\def\and{\quad{\rm and}\quad}

\def\defeq{:=}

\def\sub{\subset}

\def\and{\quad\text{and}\quad}

\DeclareMathOperator{\rank}{rank}

\let\lab=\label

\newtheorem{prop}{Proposition}[section]
\newtheorem{theo}[prop]{Theorem}
\newtheorem{lemm}[prop]{Lemma}
\newtheorem{coro}[prop]{Corollary}
\newtheorem{defi}[prop]{Definition}
\newtheorem{conj}[prop]{Conjecture}

\theoremstyle{definition}

\def\Ob{\cO b}

\def\loc{_{\mathrm{loc}}}

\def\Po{{\mathbb P^1}}

\def\Pf{{\mathbb P}^4}
\def\Pn{{\mathbb P}^n}
\def\PP{{\mathbb P}}

\def\AA{{\mathbb A}}

\def\bMPdg{\overline{M}_g(\PP^n,d)}

\def\fM{\mathfrak{M}}

\def\MPdg{M_g(\PP^n,d)}

\def\sta{^\ast}

\let\lab=\label

\def\lab#1{\label{#1}[{#1}]\  }

\def\lab{\label} %{{\bf index}-}%=\label

\def\beq{\begin{equation}}
\def\eeq{\end{equation}}

\def\bC{{\mathbf C}}

\def\tcY{\tilde{\cY}}

\def\tfd{{\tilde \cD}}

\def\cpx{{\tilde \cY}}

\def\tcy{{\tcY}}

\def\beq{\begin{equation}}
\def\eeq{\end{equation}}

\let\lab=\label

\begin{document}

\title{Relative Resolution and Its Applications}

\date{}
\author{Yi Hu}
\address{Department of Mathematics, University of Arizona, USA.}
\email{yhu@math.arizona.edu}

\maketitle

\begin{abstract}
We  present an introduction to the {\it derived and relative resolutions} of the moduli of stable maps.
We discuss one application and mention a few problems.
\end{abstract}

\bigskip\medskip

%\hskip 7.3cm {\footnotesize \it Regularities are all alike; every}

%\hskip 6.3cm {\footnotesize \it  singularity is singular in its own way.}

%\bigskip

\section{Introduction}

We aim to provide a {\it relative resolution} of a singular moduli, preferably by geometric method, not by any algorithm. Our primary example
is the Kontsevich  moduli space $\bMPdg$ of degree $d$ stable maps from genus $g$ curves to the projective space $\Pn$.
In this note, we will give a hopefully leisurely introduction to some recent works on the two kinds of $``$resolutions$"$ of
 $\bMPdg$ %of degree $d$, genus-$g$ stable maps  into $\Pn$
 (\cite{VZ08}, \cite{HL10}, \cite{HL11}, \cite{HL12}).

 The first is the {\it derived resolution}. % We will give its definition in the main text.
 Its primary purpose is to define the {\it reduced }GW invariants  which
properly  count the contributions of the main components of the moduli of stable maps.
 The derived resolution always exists and has a  minimal one,  unique up to isomorphism.

 The second is the {\it relative resolution}.   For a {\it smooth} stack $\fM$,
 we say that a blowup $\widetilde\fM \to \fM$ is smooth if $\widetilde\fM$ is smooth. Thus, smooth blowups
 include blowups along smooth closed centers.   If a moduli space has several components,
 it usually comes equipped with a distinguished one, called the main component.

 \begin{defi}\lab{rel-d} Let $X$ be a Deligne-Mumford stack with a main component $X'$.
 Let $\fM$ be a smooth Artin stack and
 $X \to \fM$ be a morphism. Suppose that we have a sequence of smooth blowups $\widetilde\fM \to \fM$
  such that if let $\widetilde{X}=X \times_{\fM} \widetilde\fM $, then
  the main component $\widetilde{X}'=X' \times_{\fM} \widetilde\fM $ is smooth and
  the entire stack $\widetilde{X}$ has at worst normal crossing singularities.
  %every other irreducible component of $\widetilde{X}$
  %is a  smooth DM stack  and they meet in normal crossing way.
  Then we say $\widetilde{X} \to X$
  is a  resolution of $X$ \underline{relative} to $\widetilde\fM \to \fM$.
   If in addition, {\it every irreducible
  component} of $\widetilde{X}$ is  smooth, then we say it is a {\it strong} relative resolution.
   \end{defi}

   Here, a DM stack (not necessarily irreducible) is said to have at worst normal crossing singularities  if at any singular point, it is locally
   equivalent to a union of coordinate subspaces of an affine space.
  A relative resolution is a natural extension of the ordinary resolution in that
  we allow  $\widetilde{X}$ to be reducible but  require normal crossing singularities for all irreducible components as well as for their
   intersections, which seems to be the best one can hope for.

   For the moduli space $\bMPdg$, one should interpret a relative resolution as follows. The stack  $\bMPdg$ fibers over a smooth stack
   $\fM$ such that all its obstructions lie in the fiber direction. A relative resolution is a process to remove all the obstructions in the fiber direction
   and at the mean time preserve smoothness of the base (whence the requirement   that the blowup $\widetilde\fM \to \fM$ is smooth).

   For many topological applications, lesser resolution may be sufficient. Let $\Eb$ be a  two-term perfect
 derived object over  the DM stack $X$.  Suppose that  $\widetilde\fM \to \fM$ is a blowup (not necessarily smooth)
 such that  if let $\widetilde{X}=X \times_{\fM} \widetilde\fM $, then
the pullback of $\Eb$ to $ \widetilde{X}  $ becomes {\it locally diagonalizable}. Then we say $\widetilde{X} \to X$
  is a {\it  derived resolution} for the object $\Eb$.  We let $\pi: \cC \to \bMPdg$ be the universal family and $\ff: \cC \to \Pn$ be the universal map.
  The canonical derived objects $\bR\pi_*\ff^*\sO_{\Pn}(k)$ are of central importance in the Gromov-Witten theory.

   Our main results may be summarized as follows.
   {\it   We obtained a strong relative resolution for $\overline{M}_1(\Pn,d)$  {\rm ([VZ08], [HL10])}.
 We achieved relative resolutions for $\overline{M}_2(\Pn,d)$ {\rm  ([HL12])}.  We constructed the derived resolutions for
 the  canonical  objects $\bR\pi_*\ff^*\sO_{\Pn}(k)$ over  $\bMPdg$ for all $g$ and $k>0$   {\rm ([HL11])}.}
More details are as follows.   %We will describe a relative resolution of $\bMPdg$ for $g=1$ and 2.

 When $g=1$,  the relative resolution is achieved over a canonical smooth blowup of
 the smooth Artin stack  of weighted curves (\cite{HL10}); for $g=2$,
  it is relative to a canonical smooth blowup of the smooth Artin stack of  pairs of curve and line bundle (\cite{HL12}).
 The resolution of the main component of $\overline{M}_1(\Pn,d)$ is originally due to Vakil and Zinger (\cite{VZ08})
and was applied to calculate the genus one reduce GW numbers by Zinger  (\cite{Zin09}).
Chang and Li  used the relative resolution to give an algebro-geometric proof of the Li-Zinger formula for genus one GW invariants of quintics
(\cite{CL12}).  It is expected that the relative resolution of  $\overline{M}_2(\Pn,d)$ can be applied to prove
the genus two version of the LZ  formula.   Higher genus cases are considerably harder. % (but seem manageable now).

We hope that this note will benefit researchers in the Gromov-Witten theory as well as other mathematicians wandering around.
All the works on the derived resolution, reduced invariants, and relative resolutions of stable map moduli
are jointly obtained with Jun Li
(\cite{HL10}, \cite{HL11}, \cite{HL12}).
Section \ref{applications} on the genus one LZ formula  is due to a joint work of Chang and Li (\cite{CL12}). I thank Jun Li for collaboration
and Huai-Liang Chang for spending hours in a series of
seminars explaining their works.   This note is based upon my lectures at Taiwan University and at
Hong Kong University of Science and Technology  during the summer of 2013.
I thank TIMS of Taiwan University for partial financial support and Chin-Lung Wang for his hospitality when I was there.
 I thank ICCM committee for its kind invitation to the sixth ICCM-Taiwan and for its subsequent  invitation to writing up my lecture notes.
 I thank Huai-Liang Chang and Wei-Ping Li for  the invitation to visit HKUST
in the summer of 2013 and for their warm hospitality; the partial financial support from HKUST is also  gratefully acknowledged.
Some of the research work described in this article was partially supported by NSF grant DMS
0901136.

The note is organized as follows.

\tableofcontents

\section{Derived and Relative Resolutions: an overwiew}

\subsection{A little background. }
The moduli space $\bMPdg$ in algebraic geometry was introduced by Kontsevich.
Its original primary application is to laid a foundation for the Gromov-Witten theory.
%Many moduli problems are important subjects of  studies in algebraic geometry and related fields.
%Among them, we just mention two:  the moduli space of stable algebraic curves and its generalization the moduli spaces of
%stable maps. The formal is smooth (as a Deligne-Mumford stack), the latter can contain arbitrary singularity types
%according to a theorem of Vakil (based upon Mn\"ev's universality theorem).
%Let  $\bMPdg$ be the moduli space of stable maps of degree $d$ from genus $g$ curves to the projective space $\PP^n$.

The moduli space $\bMPdg$ is a {\it moduli}
compactification of $\MPdg$, the open subset of  degree $d$ maps  from smooth projective curves of genus $g$ to the projective space $\PP^n$.
It is not a compactification in the usual sense since $\bMPdg$ may contain many components of various dimensions. The closure
of $\MPdg$, is the main component, denoted  $\bMPdg'$. It is of a central interest in the GW theory to understand the
 contribution to the GW invariants
of each individual component. For this, it is helpful to understand the structures of the components and how they interact, to certain extent.
But, when $n \ge 2$, the spaces $\MPdg$ (hence also $\bMPdg$) can contain arbitrary singularity types when $g$ varies and $0 < d < 2g-2$
according to a theorem of Vakil (based upon Mn\"ev's universality theorem).
%The moduli space $\bMPdg$ has a distinguished main component $\bMPdg'$ whose general point is a stable map with smooth domain curves. (We note here that $\bMPdg'$ needs not to be irreducible if $d \le 2g-2$.)
In addition, the main component $\bMPdg'$  needs not to be irreducible if $d \le 2g-2$.
So, it looks rather hopeless to fully grasp the geometric structures of the moduli space.

However, thinking positively, the moduli spaces $\bMPdg$ provide us a single geometric setting
 to resolve all possible singularities in algebraic geometry, possibly by a uniform geometric method.
%from the pure perspective of resolution singularities (see section \ref{pmr}),
% it is desirable to find a method to resolve all the singularities of $\bMPdg$.
Even if one is not interested in resolution of general singularities but is only interested in the GW invariant of hypersurfaces
or more generally the GW invariants of  complete intersections in $\Pn$, it is still
% On the other hand, for some theoretical and computational purposes in GW theory, it is also
 desirable to have a resolution of  the singularities of $\bMPdg$,  to certain extent\footnote{ E.g., in the sense of Definition \ref{rel-d}, or some weaker ones to be explained in what follows.}.  We will explain and hopefully convince the reader that a {\it partial}
 relative resolution of  $\bMPdg$, in a suitable sense, can be applied to study the GW  numbers of a quintic threefold
and more generally to study the GW invariants of a hypersurface of $\Pn$. % to the moduli $\overline{M}_g(\Pn,d)$.

 \subsection{Where are the singularites? }  So, to get a $``$resolution$"$ of singularities of the moduli space $\bMPdg$,
 first we need to understand where its singularities are.

  A point of the moduli space  $\bMPdg$ is a degree $d$ stable morphism $u: C \to \PP^n$,
 where $C$ is a projective  curve of genus $g$ with at worst nodal singularities.
Its stability means that its automorphism group is finite.
We will often abbreviate a stable map $u: C \to \PP^n$ as $[u, C]$.
 Take any point $[u_0: C_0]$ in  $\bMPdg$.
Fix  homogeneous coordinates $[x_0,\cdots,x_n]$ of $\PP^n$ with
$x_i \in H^0(\sO_{\PP^n} (1))$. Then we can write $u_0=[s_0,\cdots, s_n]$ where
$s_i=u^* x_i \in H^0(u_0^*\sO_{\PP^n} (1))$. Thus, we see that a deformation of  $[u_0, C_0]$ is determined by a combined deformation
of the curve $C_0$ and the sections $\{s_0,\cdots, s_n\}$. As the deformation of the (nodal) curve $C_0$ is unobstructed,  all possible obstructions
to deform the map  $[u_0, C_0]$ come from the obstructions of deforming the sections $\{s_0,\cdots, s_n\}$.
If in a neighborhood of  $[u_0, C_0]$, the rank of $H^0(u^*\sO_{\PP^n} (1))$ remains constant, then one has a vector bundle over the neighborhood
with fibers $H^0(u^*\sO_{\PP^n} (1))$. After shrinking the neighborhood, we may assume that this vector bundle is trivial; after a trivialization,
one sees that there is no obstruction to choose a set of sections $\{S_0,\cdots, S_n\} $ of this trivialized bundle,  extending the vectors
$\{s_0,\cdots, s_n\} $ in $H^0(u_0^*\sO_{\PP^n} (1))$ to its nearby fibers $H^0(u^*\sO_{\PP^n} (1))$. Or else,
if the rank of $H^0(u_0^*\sO_{\PP^n} (1))$ is higher than its  its nearby general  fibers $H^0(u^*\sO_{\PP^n} (1))$,
then some sections of $\{s_0,\cdots, s_n\} $ may experience obstructions when we try to extend them to nearby   fibers. If this happens,
we encounter a singularity of $\bMPdg$ at $[u_0, C_0]$.

{\it Informally, we say that a stable map has two directions: the curve direction and the section direction;
its curve direction is unobstructed but its section direction may experience  obstructions.}

Now, to formulate the above informal discussions in formal terms,
we need the universal curve $\pi: \cC \lra \bMPdg$ over $\bMPdg$
and the universal map $\ff: \cC \lra \PP^n$. For a point $[u, C] \in \bMPdg$, the fiber $\pi^{-1}( [u, C] ) \sub \cC$ is just the curve $C$;
the map $f$ restricted to the fiber $\pi^{-1}( [u, C] )=C$ is just the map $u: C \lra \PP^n$.  To capture the variation of vector space
 $H^0(u^*\sO_{\PP^n} (1))$, we have the direct image sheaf  $\pi_*\ff^*\sO_{\PP^n}(1) $.  In diagrams, this is
%\begin{equation}
$$
\xymatrix{
& & \ff^*\sO_{\PP^n}(1)  \ar[dl] &  \sO_{\PP^n}(1)  \ar[dl]  \\
\pi_*\ff^*\sO_{\PP^n}(1) \ar[dr] & \cC \ar[r]^\ff \ar[d]_{\pi } & \PP^n \\
& \bMPdg   }
$$
%\end{equation}

By the above discussions, we have the following important observation.
 The moduli space $\bMPdg$ is smooth at a point $[u, C] $ if and only if the sheaf  $\pi_*\ff^*\sO_{\PP^n}(1) $ is locally free
at $[u, C]$. Or, equivalently, the moduli space $\bMPdg$ is singular at a point $[u, C] $
if and only if the sheaf  $\pi_*\ff^*\sO_{\PP^n}(1) $ is not locally free
at $[u, C] $. (By this characterization, all smooth stable maps must lie in the main component of $\bMPdg$. The remainder components contain no smooth points.) Thus, intuitively, one would expect that {\it making the sheaf  $\pi_*\ff^*\sO_{\PP^n}(1) $  locally free will have
consequence on the resolution of the moduli space $\bMPdg$.}  We will elaborate this the later sections.

Now, what about the other direct image sheaves $\pi_*\ff^*\sO_{\PP^n}(k) $, $k >1$?

\subsection{The canonical derived objects $\bR\pi_*\ff^*\sO_{\PP^n}(k) $ }
Well, if we consider a hypersurface $X$ of degree $k$ in $\PP^n$ (e.g., a quintic in $\Pf$), we may assume that
it is given by the zeros of a section $s \in H^0(\sO_{\Pn }(k) )$. Then one checks that at least set-theoretically,   we have
that the moduli space $\overline{M}_g(X,d) \subset \bMPdg$ is given by the zeros of the section $\sigma=\pi_*\ff^* s$
of the sheaf $\pi_*\ff^*\sO_{\PP^n}(k) $. Thus, in an ideal situation such as if $\pi_*\ff^*\sO_{\PP^n}(k) $ were locally free and
$\bMPdg$ were of pure and expected dimension,  then $\pi_*\ff^*\sO_{\PP^n}(k) $
 would possess the Euler class
and the fundamental class of $\overline{M}_g(X,d) $ would be the intersection of the Euler class with the fundamental class of $\bMPdg$.
This motivates that we should also  make $\pi_*\ff^*\sO_{\PP^n}(k) $ locally free so that we can define Euler class and obtain invariants when restricted to
the main component of $\bMPdg$.

Thus, for either  purpose of desingularizing the moduli space $\bMPdg$ and
making applications to the GW theory of hypersurfaces in $\Pn$,
one would like to resolve the sheaves $\pi_*\ff^*\sO_{\PP^n}(k) $, $k \ge 1$. However, we should point out here that
as it turns out, it is not very helpful to just
consider the direct image sheaf $\pi_*\ff^*\sO_{\PP^n}(k) $,  it is more natural and mathematically correct to also include the higher
direct image sheaf $R^1 \pi_*\ff^*\sO_{\PP^n}(k) $ and  consider the derived object  $\bR\pi_*\ff^*\sO_{\PP^n}(k) $.

To stress the points, one should bear in mind that when $k=1$ the structures of the object
$\bR \pi_*\ff^*\sO_{\PP^n}(1) $  implies structural properties of the underlying moduli space $\bMPdg$; when $k>1$, the structures of the object
$\bR \pi_*\ff^*\sO_{\PP^n}(k) $  will have applications to the GW invariants of the hypersurfaces of degree $k$ in $\PP^n$.
Our treatment for $\bR \pi_*\ff^*\sO_{\PP^n}(k) $, $k \ge1$,  is uniform.

\subsection{The quintic Calabi-Yau threefolds}
 We now spend a subsection to discuss the important case of quintic CY threefolds which is related to the object $\bR\pi_*\ff^*\sO_{\Pf}(5)$.
 Let $Q \subset \Pn$ be a smooth quintic threefold. The (virtual) numbers of curves in Q are
central topics in the GW theory for the last two decades.  To this date, only the cases of  rational and elliptic curves are mathematically known.
The case of higher genus is still out of reach.

Consider the moduli space $\overline{M} _g(Q,d)$ of degree $d$ stable maps from genus $g$ curves into the quintic $Q$. Its expected dimension is zero and possesses a virtual fundamental cycle $[\overline{M} _g(Q,d)]^{\rm vir}$ which is a zero cycle with rational coefficients. Its degree is a rational number, denoted $N_{g,d}$. This provides a virtual count of genus $g$, degree $d$ curves in $Q$. Let us know examine the meanings of $N_{g,d}$
via the structures of the moduli space $\overline{M}_g(Q, d)$. The space $\overline{M}_g(Q, d)$ consists of stable maps classified as in
 Figure 1.  The domain of every stable map is a genus $g$ curve. But in its image,
 one may only see a genus $h$ curve if a genus $(g-h)$ subcurve is contracted
(the dotted curve in the picture is contracted); such a stable map is  a general point in its corresponding component $\overline{M}_g(Q, d)_h$ of  $\overline{M}_g(Q, d)$. Intuitively, this indicates that the virtual number $N_{g,  d}$ should contain  the contribution of the component of
$\overline{M}_g(Q, d)_h$; this contribution should come with two parts: one  for $\overline{M}_h (Q, d)$,
i.e., the genus  $h$ curves that we see in its images, which we denote by
 $N_{h, d}'$; the other comes from the contracted curves of genus $g-h$, and should be an invariant on
 the moduli space $\overline{M}_{g-h}$ of stable curves of genus $(g-h)$,
 which we denote by $c_{g-h}$.

 \vskip 4.5cm

\begin{picture}(3, 15)
\put(20,2){ \psfig{figure=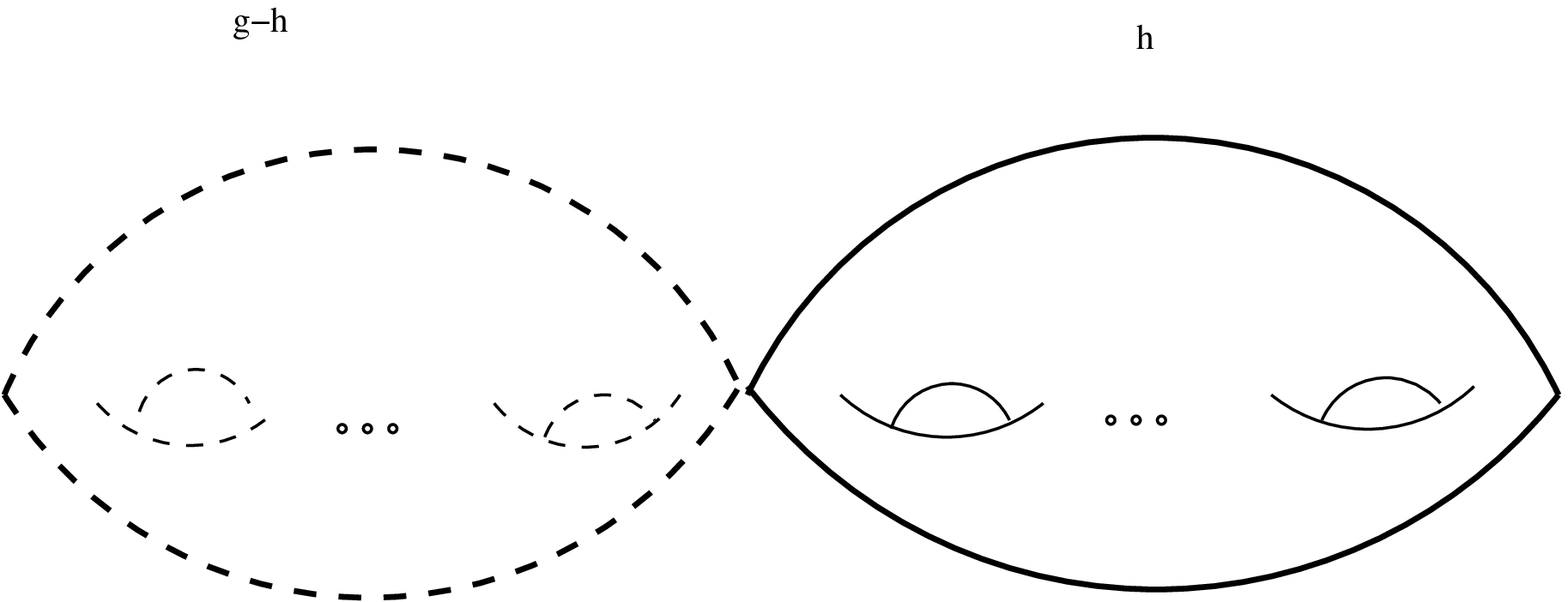, height = 4cm,width=10cm}}
\end{picture}

\vskip .6cm
\centerline{Figure 1}
\vskip .5cm

 There are other components, but conjecturally they do not contribute.  This explains the idea behind Li-Zinger's conjecture:
\beq\lab{splittingFormula}
 N_{g,  d} = N_{g,d}' + c_1 N_{g-1,d}' + \cdots + c_{g-1} N_{1,d}'  + c_g N_{0,d}'. \eeq
What is $N_{h, d}' $? This should be the contribution of the main component $\overline{M}_h (Q, d)'$ of $\overline{M}_h (Q, d)$ whose general points
are stable maps with smooth domain curves.
How do we rigorously define $N_{g,d}' $?  We will explain this in the next couple of sections.

 \subsection{The derived resolution}\lab{a-d}
 To rigorously define $N_{g, d}' $, we need to make the sheaf $\pi_*\ff^*\sO_{\Pf}(5)$ locally free so that it can provide an Euler class that we can use.
 As pointed out earlier, it is more correct to consider the object $\bR\pi_*\ff^*\sO_{\Pf}(5)$ rather than just the sheaf $\pi_*\ff^*\sO_{\Pf}(5)$ alone. We can do this in more general setting. So now we examine the properties of the object $\bR\pi_*\ff^*\sO_{\PP^n}(k)$.
 First, it is a particularly nice kind of object: it is a two term perfect object, meaning locally, it can be presented as a two term complex $[E \stackrel{\varphi }\lra F ] $ of locally free sheaves $E$ and $F$.  (Even though, for $\bR\pi_*\ff^*\sO_{\PP^n}(k) $,
  a global presentation can be found, we still prefer its local presentations, as in practice,
   local presentations come more naturally, and more importantly, possess geometric meanings.)
   Recall that one of our naive aim is to make the sheaf $\pi_*\ff^*\sO_{\PP^n}(k)$ locally free.  As it turns out, the more functorial way to achieve this
   is to {\it locally diagonalize}  the object  $\bR\pi_*\ff^*\sO_{\PP^n}(k)$ from which the local freeness of  $\pi_*\ff^*\sO_{\PP^n}(k)$ comes as an
   immediate  consequence. So, we explain now the very useful notion of {\it locally diagonalizable} derived object.

For any scheme $X$,  a  homomorphism $\varphi: \sO_X^{\oplus p} \lra \sO_X^{\oplus
q}$ is said to be diagonalizable if we have direct sum
decompositions by trivial sheaves \beq\lab{decompositionE-F}
\sO_X^{\oplus p}=G_0 \oplus \bigoplus_{i=1}^l G_i \and
\sO_X^{\oplus q}=H_0 \oplus \bigoplus_{i=1}^l H_i \eeq with
$\varphi (G_i) \subset H_i$ for all $i$ such that
\begin{enumerate} \item
$\varphi|_{G_0}=0$; \item  for  every $1\le i \le l$,
$\varphi|_{G_i}$ equals to $p_i I_i$ for some $0 \ne p_i \in
\Gamma(\sO_X)$  where $I_i: G_i \to H_i$ is an isomorphism; \item
{ $\langle p_i \rangle   \supsetneqq \langle p_{i+1} \rangle$}.
%{ The pair $p_i$ and $p_{i+1}$ are not differed by a unit in the
%structure sheaf} and $p_i | p_{i+1}$ $(1 \le i \le l-1)$.
\end{enumerate}
   A homomorphism $\varphi: E \to F$ between locally
free sheaves of a scheme $X$ is said to be locally diagonalizable if there
are trivializations of $E$ and $F$ over some open covering of $X$
such that  $\varphi: E \to F$ is diagonalizable over every open
subset.

Now if $\Eb$ is a two-term perfect derived object over a DM stack $X$, using any local presentation of  $[E \stackrel{\varphi }\lra F ] $
of $\Eb$, we can introduce the notion of locally  diagonalizable derived
object.  One can prove that the notion does not depend on the local presentation
(\cite{HL11}).

Not every homomorphism (or two-term perfect derived object)  can be locally diagonalized.
But if it can be, then there are some good implications.
For example, suppose that a homomorphism $\varphi: E \to F$ (or an object $\Eb$) is  locally
diagonalizable, then for every irreducible
component $X'$ of $X$ with the reduced scheme structure, $\ker(\varphi|_{X'})$ ($\sH^0(\Eb)$) is locally free. Also, if $f: Y \to X$ is a
morphism, then $f^* \varphi$ ($f^* \Eb$) is also locally diagonalizable, that is,
$``$locally diagonalizable$"$ has base change property.
Note here that the rank
of $\ker \varphi$ depends on the properties of the functions
$p_i$, hence may not be constant over $X$. { Further, when
$X'$ is not reduced, $\ker (\varphi|_{X'})$ needs not to be
locally free. These technical issues lead us to use the
$``$integral$"$ assumption whenever  we want to
produce a locally free sheaf.

 For any $0 \le i \le r:={\rm max}(\rank E,\rank F)$,
 we  let $\sI_{\varphi, i}$ be the $i$-th determinantal ideal sheaf of $\varphi: E \to F$.
  That is, $\sI_{\varphi, i}$ is the ideal sheaf of the zero scheme of $\wedge^{i+1} E \to
  \wedge^{i+1}F$. Again, using local presentations of the object $\Eb$, one can introduce
  the $i$-th determinantal ideal sheaf $\sI_{\Eb, i}$ of $\Eb$ and show that  it is independent of the choice of the
  presentations.  We  proved the following.

\begin{theo}{\rm (\cite{HL11})}\lab{derivedR}  Suppose $\Eb$ is a two-term perfect derived object
over a DM stack $X$.  If we let $b: \tilde X  \lra X$ be the blowup of $X$
  along the ideal sheaves  $\sI_{\Eb, 0}, \cdots, \sI_{\Eb, r}$ and so on, then the pullback $b^* \Eb$
  %$b^*\varphi: b^*E \lra  b^*F$
  becomes locally diagonalizable over $\tilde X$.  In particular, $\sH^0(b^*\Eb)$
is locally free provided $X$ is integral.
  %In particular, $\ker b^*\varphi:$ is  locally free provied $X$ is integral.
  Further, $\tilde X$ is minimal in the sense that if there is another dominating
 morphism $g:  Z \to X$ such that $g^*\Eb$ is locally diagonalizable, then $g$ factors uniquely as $Z \to \tilde X \to X$.  %Further, the blowup $b: \tilde X  \lra X$ is minimal in the following sense:
  %if $f: Z \lra X$ is another dominant morphism such that   $f^*\varphi: f^*E \lra  f^*F$ is locally diagonalizable,
  %then $f$ factors as $ Z \lra \tilde X \lra X$.
   \end{theo}

  In fact, we can consider any perfect derived object $\sEb$ in the
 bounded derived category  ${\rm D}^b(M)$ of an
integral DM stack $M$ with  cohomologies concentrated in the non-negative
places. Using the above theorem, we can show

\begin{theo}\lab{Dresolution}  {\rm (\cite{HL11})}  Let $\sEb$ be any perfect derived object over an integral DM
stack $X$.  Assume that $\sEb$ can  locally be represented
by a complex of locally free sheaves of finite length supported
only in non-negative degrees. %$\sH^i(\sEb)=0$ for $i<0$.
Then there is another
integral DM stack $\widetilde{X}$ and a surjective birational
morphism $b: \widetilde{X} \to X$ such that  $\sH^0( b\sta \sEb)$ is
locally free.
\end{theo}

The above allows us to define the Euler class of such a derived object.
For this, we suppose further that the cohomologies in positive places
$\sH^{i>0}(\sEb)$ are all torsion sheaves  over $X$.
Then, we define the Euler class ${ \ee}(\sEb)$ in the Chow group $A_* X$  of cycles on
$X$ by, \beq { \ee}(\sEb) := b_*(c_r(\sH^0(b\sta \sEb)) \cdot
[\widetilde{X}]),\eeq  where $r=\rank  \sH^0( b\sta \sEb)$.
One can show the Euler class ${ \ee}(\sEb)$ is independent of the choice of the
resolutions $b: \widetilde{X} \to X$.

\subsection{The reduced GW numbers of quintics and LZ conjecture}
The above, when applied to the derived object
$\bR\pi_* \ff^* \sO_{\PP^4}(5)$ restricted to
the primary component $\overline{\fM}_g(\PP^4, d)'$
 of the moduli stack
$\overline{\fM}_g(\PP^4,d)$, enables us to construct the modular Euler class when $d > 2g-2$.
Here, the general points of $\overline{\fM}_g(\PP^4,d)'$ are
maps with smooth domains, and
%general smooth domain curves such that the maps are defined by sections of
%general divisors on the curve.
$\overline{\fM}_g(\PP^4,d)'$ is irreducible and of the expected dimension when
$d >2g-2$.   In this case, letting $\pi'$ and $f'$ be the restrictions of the projection of the universal family and universal map
to $\overline{\fM}_g(\PP^4,d)'$, respectively,  then $R^1 \pi'_* f^{\prime\ast} \sO_{\PP^4}(5)$
is a torsion sheaf over $\overline{\fM}_g(\PP^4,d)'$. Hence,   for $d> 2g-2$,
we can define the modular Euler class of $\bR\pi'_* \ff^* \sO_{\PP^4}(5)$ over $\overline{\fM}_g(\PP^4,d)'$ to be
\beq { \ee}(\bR\pi'_*f^{\prime\ast}
\sO_{\PP^4}(5)) \in A_* (\overline{\fM}_g(\PP^4,d)');
\eeq
for any smooth Calabi-Yau manifold $Q$ in $\PP^4$, we define
\beq N_{g,d}'(Q) = \deg \ee(\bR\pi'_* f^{\prime\ast} \sO_{\PP^4}(5)).
\eeq

Using these numbers $N_{g,d}'(Q)$, in \cite{HL11}, we  rigorously formulate  the conjecture of Li and Zinger:

\begin{conj}\lab{conj}  Let $N_{g,d}(Q)$ be the genus $g$, degree $d$ GW-invariants of the smooth
quintic $Q\sub\PP^4$.%assume that the primary  part $\overline{\fM}_g(\PP^4,d)'$ is not empty and every of its components
%has the pure expected dimension. Then
Then there are universal constants $c_h$
such that for all  $d > 2g-2$,
$$N_{g,d}(Q)= \sum_{0 \le h \le g} c_h N_{h,d}'(Q).$$
\end{conj}

The $g=0$ case is trivial: $N_{0,d}=N_{0,d}'$.
The $g=1$ version, $N_{1,d}=N_{1,d}'+\frac{1}{12} N_{0,d}'$,
 is proved by Li and Zinger (and later by Chang-Li, using algebro-geometric method). The $g \ge 2$ cases are  open.

%\subsection{Setting up the stage for relative resolution}

\subsection{Relative resolutions: a general program}\lab{rrForAll}
 The derived resolutions, i.e., the blowup $b: \tilde X  \lra X$  in    Theorem \ref{derivedR}
(or Theorem \ref{Dresolution}), are unfortunately very singular in general. This is so even when we start with a smooth scheme $X$. For various purpose, a  resolution in the usual sense,  which, at the mean time,  can also serve as a derived resolution,
 is desirable. This is, of course, much hard to achieve in general.

For the case of $\bMPdg$, we hope to achieve a relative resolution.
Recall that a stable map has two directions to deform, one is  curve, the other is  collection of sections. Intuitively,
one could image that the derived resolution removes the obstructions of the section directions but at the cost of destroying the (original)
smoothness of the curve direction. The difficulty is to resolve the obstructions of
section direction and at the mean time preserve the smoothness of the curve direction. Below, we will make
it clear this vague point.

To find out the base of relative resolution  of $\bMPdg$,
%we first set up the stage for the relative resolution of $\bMPdg$.
we should begin with  analyzing a neighborhood $U$ of a singular point $[u,C] \in \bMPdg$. We  use $\pi_U:
\cC_U \to U$ to denote the restriction $\cC|_{U}$ of the universal family and $f_U$ the restriction $f|_{\cC_{U}}$ of the universal map. So, in diagram, we have
$$\begin{CD}
\cC_U @>{f_U}>> \PP^n \\
@V{\pi_U}VV  \\
U.
\end{CD}
$$
Suppose we plan to resolve $\bR\pi_*\ff^* \sO_{\PP^n}(k)$ (just remember that the case $k=1$ ties closely to the local structures of
the underlying moduli). First, we represent $\sO_{\PP^n}(k)$ as $\sO_{\PP^n}(H)$ by picking a general hypersurface of degree $k$ meeting
the image curve $f(C)$ (which is of degree $d$) at $m:=dk$ many {\it isolated} points. If we perturb the map $[u,C]$ slightly, the image curve still meets $H$ at isolated $m$ points. Putting all these neighboring maps together,
this implies that by replacing $U$ by a smaller neighborhood  if necessary,
we may assume that we have  $f_U^* \sO_{\Pn} (H)= \sO_{\cC_U} (D)$ with $D=D_1 + \cdots +D_m$ where
each $D_i$ is a section of the family $\cC_U \to U$ and they are  disjoint.
If $C'$ is a component of $C$ such that the degree of $u|_{C'}$ is $d'$, them there should be $d'k$ many sections passing through
the curve $C' \sub \cC_U$; if $C'$ gets contracted by the map $f$, then $C'$ meets none of the sections.

Observe that the choice of the general hypersurface $H$ determines a morphism
by assigning the map $[u,C]$ to the pair $(C, \delta_1+\cdots +\delta_m)$ where $\delta_i=C \cap D_i$.  For this purpose,
we introduce
the Artin stack $\fM^{\rm div}_g$ of pairs $(C,\delta_1+\cdots +\delta_m)$ where $C$ is a prestable curve of genus $g$  and
$\delta_1, \cdots , \delta_m$ are disjoint smooth points on $C$.
This  stack is smooth.  By the above, we have a morphism $U \to  \cV \sub \fM^{\rm div}_g$
where $\cV$ is a chart of $\fM^{\rm div}_g$. (The morphism is not canonical because it depends on the choice of $H$).
 The geometry of the pairs is easier to get hold of than that of the stable map. This leads to the idea that we should work over
the chart $\cV$ and then pull back the results to the neighborhood $U$. For this, the Artin stack admits a universal curve $\cC_\cV$
and a universal divisor $\cD=\cD_1+ \cdots + \cD_m$ so that we have a fiber square
$$\begin{CD} \cC_U @>{\alpha}>>  \cC_\cV & \supset \cD \\
@V{\pi_U}VV  @V{\rho_\cV}VV \\
U @>>> \cV
\end{CD}
$$
and $D_i=\alpha^* \cD_i$ and $D=\alpha^* \cD$. This allows us to study the object $\bR(\pi_U)_* \sO_{\cC_U}(D)$ via the object
$\bR(\rho_\cV)_* \sO_{\cC_\cV}(\cD)$.  That is, a derived resolution of $\bR(\rho_\cV)_* \sO_{\cC_\cV}(\cD)$ will induce a derived resolution of
 $\bR\pi_* \sO_{\cC_U}(D)$ because of the base change property.
 The task now is: smoothly blow up $\cV$  so that the object $\bR(\rho_\cV)_* \sO_{\cC_\cV}(\cD)$ becomes
 locally diagonalizable. This is what we mean earlier when we mention that to get a true resolution we should
 {\it  resolve the obstructions of the section directions {\rm (}i.e., to locally diagonalize $\bR(\rho_\cV)_* \sO_{\cC_\cV}(\cD)${\rm )}
 and preserve the smoothness of curve directions {\rm (}i.e., smoothly blowup the chart $\cV${\rm )}.}

 Unfortunately, we have no obvious
global morphisms from $\bMPdg$ to $\fM^{\rm div}_g$. To overcome to this difficulty,
one could in principle  work over local charts $\cV$ of $\fM^{\rm div}_g$ and then (cumbersomely) argue how to patch the data over the charts. Instead, we choose to introduce another Artin stack so that it admits a global morphsim from $\bMPdg$ and not too much
technical information is lost.

 This substituting stack is the Artin stack
$\fP_g$ of pairs $(C,L)$ of a prestable curve of genus $g$ and a line bundle over $C$. This again is a smooth stack. It comes with a universal curve
$$\rho_{\fP}: \cC_{\fP} \lra \fP_g$$ and a universal line bundle $\sL$. The relevant derived object here is $\bR(\rho_{\fP})_* \sL$.
We have a (global) morphism (depending on $k$)
$$\varsigma_k: \bMPdg \lra \fP_g, \;\;\; [u, C] \to (C, u^*\sO_{\PP^n}(k)).$$
Our earlier local morphism clearly factors through
$$
\xymatrix{
U  \ar[d] \ar[ddrr]   |\hole  \ar[rr]  &&  \bMPdg \ar@{.>}[d]  \ar@/^25pt/[dd]^{\varsigma_k}\\
\cV  \ar[rr]  && \fM_g^{\rm div} \ar[d]\\
 && \fP_g
}
$$
%$${\varsigma_k}_{U}: U (\sub \bMPdg) \lra \cV (\sub \fM^{\rm div}_g) \lra \fP_g.$$
The idea now is that we can work over the chart $\cV$ (which we found easy to cope with since it conveniently
parameterizes a pair of curve and a divisor).  Then, according the local smooth blowup of $\cV$ we find,
we may translate it to the stack $\fP_g$ to obtain the corresponding global smooth blowup, as desired.
The goal is that we need to smoothly blow up $\cV$  and hence also $\fP_g$ so that  $\bR(\rho_{\fP})_* \sL$,
 hence  also $\bR\pi_*\ff^* \sO_{\Pn}(k)$,
becomes locally diagonalizable. the author excepts the following.

\begin{conj}\lab{d-r}   Fix $g>0$ and $k>0$.
There is a canonical sequence of smooth blowups $\tilde\fP_g^k \to  \fP_g$ of $\fP_g$  such that
if we form the cartesian diagram
$$
\begin{CD}
\widetilde{M}_g^k(\Pn,d)= \tilde\fP_g^k \times_{\varsigma_k, \fP_g} \overline{M}_g(\Pn,d) @>>> \tilde\fP_g^k \\
 @VVV @VVV \\
 \bMPdg @>>> \fP_d,
 \end{CD}
 $$
then the derived object $\bR\pi_*\ff^* \sO_{\Pn}(k)$, when pulling back to $\widetilde{M}_g^k(\Pn,d)$, becomes locally diagonalizable.
\end{conj}

Such a {\it relative derived resolution} should be sufficient to apply to the GW theory\footnote{Indeed,
an easier  partial relative derived resolution is sufficient. This will be explained in  a forthcoming paper.} of hypersurfaces in $\Pn$.
Because $\tilde\fP_g^k$ is smooth, we have a perfect relative obstruction theory $\widetilde{M}_g^k(\Pn,d)/\tilde\fP_g^k$
by pulling back a perfect relative obstruction theory $\bMPdg/\fP_g$ (which is known to exist).
We may apply this to $R\pi_*\ff^* \sO_{\Pf}(5)$ to study the GW invariants of quintic threefolds.

To get a relative resolution of $\bMPdg$, we specialize the above to the case of $k=1$.

\begin{conj}\lab{r-r}  Fix $g>0$.  There is a canonical sequence of smooth blowups $\tilde\fP_g \to  \fP_g$ of $\fP_g$
 such that if we form the cartesian diagram
$$
\begin{CD}
\widetilde{M}_g(\Pn,d)= \tilde\fP_g \times_{\fP_g} \overline{M}_g(\Pn,d) @>>> \tilde\fP_g \\
 @VVV @VVV \\
 \bMPdg @>>> \fP_d
 \end{CD}
 $$
 (we omit $\varsigma_1$ from the fiber product),
 then the derived object $\bR\pi_*\ff^* \sO_{\Pn}(1)$, when pulling back to $\widetilde{M}_g^k(\Pn,d)$, becomes locally diagonalizable.
 Further, assume that $d > 2g-2$, then we have
\begin{enumerate}
\item  the main component of $\widetilde{M}_g(\Pn,d)$ is smooth;
\item  the entire stack $\widetilde{M}_g(\Pn,d)$ has at worst normal crossing singularities.
\end{enumerate}
\end{conj}

 The above is a resolution of $\overline{M}_g(\Pn,d)$ {\it relative to} $\fP_g$.
 This sort of {\it relative resolution} problem may be characterized as
$``${\it desingularization by removing relative obstructions, while preserving the smoothness of the base.}$"$
The author expects this to be true possibly after substituting $\fP_g$ by another moduli stack, or re-interpret the smoothness.

We can carry out this program for the cases of genus one and two. That is, Conjectures \ref{d-r} and \ref{r-r} hold in this two cases.
In fact, our statements are slightly stronger  than stated in the conjecture.

\subsection{Relative resolution: the genus one case}
Ideally,  for technical convenience as we discussed earlier,  we would prefer to describe a relative resolution of
$\bMPdg$ with $\fM^{\rm div}_g$ as the base.  But, we do not have a global morphism to achieve so.
Thus we are forced to compromise to work  over the stack  $\fP_g$ downstairs.
However, the stack $\fP_g$ is not that convenient to work with.  Therefore,
whenever possible, we should replace $\fP_g$  by another stack that is easier to work with.
For genus one, this substituting stack is the smooth Artin stack $\fM_1^{\rm wt}$ of weighted curves.
In general, the stack $\fM_g^{\rm wt}$ consists of pairs of pre-stable curve of genus $g$ and a nonnegative  integer (called weight). Forgetting the weight, we obtain a morphism from $\fM_g^{\rm wt}$ to the smooth Artin stack $\fM_g$  of pre-stable curves.  This morphism is \'etale, showing the smoothness of $\fM_g^{\rm wt}$.  There is a canonical morphism $\fP_g \to \fM_g^{\rm wt}$ by assigning  a pair $(C,L) \in \fP_g$ to
$(C, c_1(L)) \in \fM_g^{\rm wt}$.  Thus,  for every $k>0$, we have
$$\bMPdg \to  \fP_g \to\fM_g^{\rm wt},\;\; [u,C] \to (C, u^*\sO_{\Pn}(k)) \to (C, c_1(u^*\sO_{\Pn}(k))). $$
For genus one, the resolution can be achieved relative to  the base $\fM_1^{\rm wt}$,

In \cite{VZ08},  Vakil and Zinger  provided a natural resolution of the main component of $\overline{M}_1(\Pn,d)$.
We follow exactly their blowing up procedure,  but instead  blow up  $\fM_1^{\rm wt}$, and then  treat the entire moduli space
by taking fiber product. In this case,
we will only use the morphism
$$\varsigma_1: \overline{M}_1(\Pn,d) \lra \fM_1^{\rm wt},\;\; [u,C] \to  (C, c_1(u^*\sO_{\Pn}(1))). $$
For any $i > 0$, let $\Theta_i \subset \fM_1^{\rm wt}$ be the close subset such its general points are of the form $(C,w)$
such that $C$ can be obtained from a smooth elliptic curve $E$ by attaching (connected) trees of rational curves at $i$ many disjoint
points of $E$ and the weight $w$ restrited to $E$ is zero.  Every $\Theta_i $ is a smooth closed substack of codimension $i$.

\begin{theo}{\rm ([VZ08], [HL10])} \lab{thm:m1}
We let $\widetilde\fM_1^{\rm wt} \lra \fM_1^{\rm wt}$ be the successive blowup of $\fM_1^{\rm wt}$
 along  the smooth closed substacks $\Theta_2, \Theta_3, \cdots,$ and so on.
If we let $$\widetilde{M}_1 (\Pn,d) = \widetilde\fM_1^{\rm wt} \times_{ \fM_1^{\rm wt}}
\overline{M}_1(\Pn,d), $$ then we have
\begin{enumerate}
%\item  the main component of $\widetilde{M}_1 (\Pn,d)$ is smooth;
\item  every (including the main) irreducible component of $\widetilde{M}_1 (\Pn,d)$ is smooth;
\item  the entire DM stack $\widetilde{M}_1 (\Pn,d) $ has at worse normal crossing singularities; further
\item the derived object $\bR\pi_*\ff^*\sO_{\Pn}(k)$, upon pulling back to $\widetilde{M}_1 (\Pn,d) $, becomes locally diagonalizable for all
$k >0$; in particular
\item for any irreducible component $N$ of $\widetilde{M}_1(\Pn,d)$,
let $(\pi_N, \ff_N)$ be the pullback of the universal family $(\pi, \ff)$ to $N$,
then the direct image sheaf
$(\pi_{N})_*\ff_N^* \sO_{\Pn}(k)$ is locally free for all $k >0$.
\end{enumerate}
\end{theo}

\subsection{Relative resolutions: the genus two case}
The genus two case is already substantially more complicated than the case of genus one.
For one thing, unlike $g=1$, for $g=2$, there are more than one rounds of canonical sequence of blowups,
and the case $k=1$ requires its own  sequence of blowups;
for another, the blowup centers are complicated to describe. Because it requires substantial preparation to describe them,
%of lacking space and time,
we will  have to be vague on the blowing up centers here. But, we will give a good hint near the end of the next section.

\begin{theo}{\rm ([HL12])} \lab{thm:m2}
There are \underline{three} rounds of  sequences of successive \underline{smooth} blowups\footnote{We mention here that the first two rounds
of sequences of smooth blowups  are {\it topological} and the third round is {\it geometric}. We will be more specific
about them in the next section. We believe that the sequences of topological blowups extend to all genera
and are sufficient to study the GW numbers of quintics. This will appear in a forthcoming paper.}
 $\widetilde\fM_2^{\rm wt} \lra \fM_2^{\rm wt}$  of $\fM_2^{\rm wt}$
such that
if we let $$\widetilde{M}_2^k (\Pn,d) = \widetilde\fM_2^{\rm wt} \times_{\varsigma_k, \fM_2^{\rm wt}} \overline{M}_2(\Pn,d), $$ then we have
\begin{enumerate}
\item  the derived object $\bR \pi_*\ff^*\sO_{\Pn}(k)$, upon pulling back to $\widetilde{M}_2^k (\Pn,d) $, becomes locally diagonalizable for all
$k \ge 2$; in particular
\item for any irreducible component $N$ of $\widetilde{M}_2^k(\Pn,d)$,
let $(\pi_N, \ff_N)$ be the pullback of the universal family $(\pi, \ff)$ to $N$,
then the direct image sheaf
$(\pi_{N})_*\ff_N^* \sO_{\Pn}(k)$ is locally free for all $k \ge 2$.
\end{enumerate}
\end{theo}

Unlike $g=1$, some of the smooth blowups for $g=2$ are not blowups along smooth closed centers.

The object $\bR \pi_*\ff^*\sO_{\Pn}(1)$ underpins the local structure of the underlying moduli space and requires one more sequence of blowups to resolve.

\begin{theo}{\rm ([HL12])} \lab{thm:m3}  There is a canonical sequence of smooth blowups of
$\widetilde\fP_2 \to \fP_2 \times_{\fM_2^{\rm wt} } \widetilde\fM_2^{\rm wt} $  such that
if we let  $$\widetilde{M}_2 (\Pn,d) = \widetilde\fP_2 \times_{ \fP_2}
\overline{M}_1(\Pn,d), $$  and suppose $d >2$, then we have  \begin{enumerate}
\item  the main  component of $\widetilde{M}_2 (\Pn,d)$ is smooth;
\item  the entire DM stack $\widetilde{M}_2 (\Pn,d) $ has at worse normal crossing singularities; further
\item the derived object $\bR \pi_*\ff^*\sO_{\Pn}(k)$, upon pulling back to $\widetilde{M}_2 (\Pn,d) $, becomes locally diagonalizable for all
$k \ge 1$; in particular
\item for any irreducible component $N$ of $\widetilde{M}_2(\Pn,d)$,
let $(\pi_N, \ff_N)$ be the pullback of the universal family $(\pi, \ff)$ to $N$,
then the direct image sheaf
$(\pi_{N})_*\ff_N^* \sO_{\Pn}(k)$ is locally free for all $k \ge 1$.
\end{enumerate}
\end{theo}

 \section{Why  does Relative Resolution Work? }\lab{details}
In this section, we provide some details to explain  how we could potentially derive the modular resolutions for all $g$
and why the cases $g=1$ and 2 work.

We remind the reader that our goal is to smoothly blow up $\fP_g$ so that the derived object
$\bR(\rho_{\fP})_*\sL$ becomes locally diagonalizable. As in Section \ref{rrForAll}, locally, we can factorize the canonical morphism
$\bMPdg \to \fP_g$ as
$$U \lra \cV \lra \fP_g$$
where $U$ is an open neighborhood of a point $[u,C] \in \bMPdg$ and $\cV$ is a chart over $\fM_g^{\rm div}$.
Any blowup of $\fP_g$ induces
a blowup of $\cV$ by taking the fiber product.
By now, we get the idea that to carry out  our relative resolution program, it is important to know the local structures of the derived object
$\bR(\rho_{\cV})_*\sO_{\cC_\cV}(\cD)$  (which is related to $\bR\pi_*\ff^* \sO_{\Pn}(k)$ by Cohomology and Base Change).

So, we begin to analyze  the derived object $\bR(\rho_{\cV})_*\sO_{\cC_\cV}(\cD)$ and pin down a local presentation.
This is known to admit a local presentation of the form $[E \stackrel{\psi}{\lra} F]$ where $E$ and $F$ are locally free sheaves.
But now, we are not content with such an inexplicit form (for the derived resolution, it suffices since we blow up along determinantal ideals
whose explicit forms are not necessary), we need an explicit form
so that we can extract the geometry behind the form. For this purpose,  shrinking $\cV$
if necessary, we can choose $g$ many disjoint sections $\cA_1, \cdots, \cA_g$ of $\cC_\cV \to \cV$ such that they are disjoint from the divisor $\cD$
and for any $v \in \cV$ and any
 irreducible  subcurve $C'$ of $C_v =\rho_\cV^{-1}(v)$ of genus $h \le g$, the set
$\{\cA_1 \cap C', \cdots, \cA_g \cap C'\}$ consists of $h$ many distinct
points in general position (that is, they do not form a special divisor on the curve $C'$).
Set $\cA= \sum_i \cA_i$. Then from the short exact sequence,
$$0 \lra \sO_{\cC_\cV}(\cD) \lra \sO_{\cC_\cV}(\cD +\cA) \lra \sO_{\cA}(\cA) \lra 0$$
we have a long one
$$(\rho_{\cV})_*\sO_{\cC_\cV}(\cD) \hookrightarrow  (\rho_{\cV})_*\sO_{\cC_\cV}(\cD +\cA)  \stackrel{\psi}{\lra}  (\rho_{\cV})_*\sO_{\cA}(\cA) \twoheadrightarrow R^1(\rho_{\cV})_* \sO_{\cC_\cV}(\cD) $$
because one shows  that $R^1(\rho_{\cV})_* \sO_{\cC_\cV}(\cD+ \cA) =0$ by the assumptions on $\cA$. Hence,
$(\rho_{\cV})_*\sO_{\cC_\cV}(\cD +\cA) $ is locally free of
rank $d+1$ by Riemann-Roch. Also, it is clear that $(\rho_{\cV})_*\sO_{\cA}(\cA) $ is locally free of rank $g$ since $\cA$ are disjoint sections.
Thus, we have
$$\bR(\rho_{\cV})_*\sO_{\cC_\cV}(\cD) = [(\rho_{\cV})_*\sO_{\cC_\cV}(\cD +\cA)  \stackrel{\psi}{\lra}  (\rho_{\cV})_*\sO_{\cA}(\cA)].$$
We will call $\psi$ a structural homomorphism;  understanding this homomorphism
is a key to our program.
We can always choose frames of  the  free sheaves
$$E=(\rho_{\cV})_*\sO_{\cC_\cV}(\cD +\cA)  \and F= (\rho_{\cV})_*\sO_{\cA}(\cA)$$
so that  the homomorphism $\psi$ can be represented by  a matrix. Then our task becomes to
derive an explicit form for every entry in the matrix. Before we do it, we need some  useful lemmas to decompose  the homomorphism $\psi$.

First we choose another general section $\cB$; this will help us to isolate a canonical trivial sheaf in the kernel of $\psi$.

\begin{lemm} We have a natural splitting of locally free sheaves %(shinking $\cV$ if necessary)
$$(\rho_{\cV})_*\sO_{\cC_\cV}(\cD +\cA) = \sO_\cV \oplus (\rho_{\cV})_*\sO_{\cC_\cV}(\cD +\cA-\cB)$$ such that
$\sO_\cV \sub \ker \psi$ and $\psi$ can be naturally  decomposed  as
$$\sO_\cV \oplus (\rho_{\cV})_*\sO_{\cC_\cV}(\cD +\cA-\cB)  \stackrel{0 \oplus \varphi}{\lra}  (\rho_{\cV})_*\sO_{\cA}(\cA) $$
where $ \varphi:  (\rho_{\cV})_*\sO_{\cC_\cV}(\cD +\cA-\cB)   \to  (\rho_{\cV})_*\sO_{\cA}(\cA) $ is the restriction homomorphism.
\end{lemm}

Since $\psi=0 \oplus \varphi$, we also call $\varphi$ a structural homomorphism.
Observe now that the rank of the locally free sheaf $(\rho_{\cV})_*\sO_{\cC_\cV}(\cD +\cA-\cB)$ coincides with the the degree of $\cD=\sum_{i=1}^m \cD_i$.  Indeed, if we let
$ \varphi_j:  (\rho_{\cV})_*\sO_{\cC_\cV}(\cD_j +\cA-\cB)   \to  (\rho_{\cV})_*\sO_{\cA}(\cA) $ be the restriction homomorphism, then we have

\begin{lemm}
We have the following natural decomposition of locally free sheaves
$$(\rho_{\cV})_*\sO_{\cC_\cV}(\cD +\cA-\cB) = \bigoplus_{j=1}^m (\rho_{\cV})_*\sO_{\cC_\cV}(\cD_j+\cA-\cB) .$$ Further, the homomorphism
$\varphi$ also decomposes as
$$\varphi= \bigoplus_{j=1}^m  \varphi_j.$$
\end{lemm}

This reduces the study of the homomorphism $\psi$ to each $\varphi_i$ which can further be reduced to
$$ \varphi_{ij}:  (\rho_{\cV})_*\sO_{\cC_\cV}(\cD_j +\cA-\cB)   \to  (\rho_{\cV})_*\sO_{\cA_i}(\cA_i)\cong \cV $$
for all $1 \le i \le g$. Shrinking $\cV$ if necessary, we can assume all sheaves $ (\rho_{\cV})_*\sO_{\cC_\cV}(\cD_j +\cA-\cB)  $
are isomorphic to the trivial sheaf $\sO_\cV$. Thus, upon fixing frames,  we may regard $\varphi_{ij}$ as a function in $\Gamma(\sO_\cV)$.
We found that this function vanishes at $v \in \cV$ if and only if there are separable nodes between $\cA_i \cap \cC_v$
and $\cD_j \cap \cC_v$. Here a node $q$ of a connected curve $C$ is separable if $C \setminus q$ is disconnected;  the separable node
$q$ is said to be between two points $a, \delta \in C$ if $a$ and $ \delta$  lie in different components of $C \setminus q$.
For any  two points $a, \delta \in C$, we let $N_{[a,\delta]}$ denote the set of all separable nodes between $a$ and $\delta$.
Every node $q$ of $C$ is associated with a function $\zeta_q$ in $\Gamma(\sO_\cV)$, called a node-smoothing parameter. Geometrically,
at any point $v \in (\zeta_q=0)$, the fiber curve $\cC_v$ contains a node corresponding to $q$; at any point $v$ away from the locus
$(\zeta_q=0)$, the curve $\cC_v$ does not have a node corresponding to $q$ (the node is smoothed out).
This leads us to the following key proposition.

\begin{prop}
Upon fixing frames,  we have  $\varphi_{ij}=c_{ij} \prod_{q \in N_{[a_i, \delta_j]} }  \zeta_q$
where $c_{ij} \in \Ga (\sO_\cV^*)$. \end{prop}

Putting all together, we obtain a matrix representation of the structural homomorphism $\varphi$ by
$$\Phi = (c_{ij} \prod_{q \in N_{[ a_i,\delta_j]} } \zeta_q)_{1 \le i \le g, 1 \le j\le m}.$$

This matrix is important because it provides a local equation of $\bMPdg$ at the point $[u,C]$ when
$m=d$ (i..e, for $\pi_*\ff^*\sO_{\Pn}(1)$).
The idea is roughly that the kernel of the homomorphism $\varphi$ gives rise to sections of $\pi_*\ff^*\sO_{\Pn}(1)$
but sections of $\pi_*\ff^*\sO_{\Pn}(1)$ determine stable maps.
Note that $\Phi$ is a matrix of size $g \times d$ (for $\pi_*\ff^*\sO_{\Pn}(1)$).
We let ${\bf w}^i$ be the {\it column} vector $(w^i_1, \cdots, w^i_d), i=1, \cdots, n$ with $w^i_j  \in \AA^1$.

\begin{theo}   (\cite{HL10})  The local equation of a (small) neighborhood  $U$  of
the point $[u,C] \in \bMPdg$ is given by  the system of equations
\beq\lab{localEq}
\Phi \cdot {\bf w}^i =0,  \;\;i=1, \dots, n, \eeq
realizing it as a closed subset of the smooth space $\cV \times \AA^{dn}$.
\end{theo}

The simplest kind of matrices are diagonal ones. In such a case, equations \eqref{localEq} become elementary.
However, one can not hope in general to find a diagonal form of  $\Phi $ simply by choosing frames.
To make  $\Phi $ diagonal in general, birational base change is necessary. A good thing is that the local equation \eqref{localEq}
pulls back to provide a local equation under any base change.

\begin{theo} Let $f: \cV' \to \cV$ be any morphism,  $U' = U \times_\cV \cV'$ and
$\Phi'$ a matrix representation of the pullback homomorphism  $f^* \varphi$. Then $U'$ is defined by
the system of equations
\beq\lab{localEqPullsBack}
\Phi' \cdot {\bf w}^i =0,  \;\;i=1, \dots, n, \eeq
realizing it as a closed subset of  $\cV' \times \AA^{dn}$.
\end{theo}

Now,  {\it suppose} $f: \tilde\cV \to \cV$ is a smooth blowup of $\cV$  such that the derived object
$\bR (\rho_\cV)_* \sO_{\cC_\cV}(\cD)$ (related to $\bR\pi_*\ff^*\sO_{\Pn}(1)$) becomes diagonalizable upon pulling back to
$\tilde\cV$. Then we can find a diagonal form for its corresponding matrix $\tilde\Phi$. Write $\tilde\Phi$ as
${\rm diag}(z_1, \cdots, z_g, 0, \cdots,0)$ (here, we assume $d>g$). Then

\begin{coro} %Under the above assumptions,
The open subset   $\tilde U = U \times_\cV \tilde\cV$
is given by equations
\beq\lab{localEqNormalCrossing}
z_j \cdot {\bf w}^i_j =0,  \;\; 1 \le  j \le g, \;  i=1, \dots, n, \eeq
realizing it as a closed subset of the smooth space $\tilde\cV \times \AA^{dn}$.
\end{coro}

The main component of $\tilde U$ is found by requiring the non-vanishing of $z_j$ over its generic points,
hence it is given by  $$w^i_j =0, \;\;  1 \le  j \le g, \;  i=1, \dots, n, $$
which define a smooth closed subset of the smooth space $\tilde\cV \times \AA^{dn}$. This shows the smoothness of the main component.
Since  $f: \tilde\cV \to \cV$ is a smooth blowup of $\cV$,
we conclude that $(z_j=0)$ defines a normal crossing divisor over $\cV$.
Any other irreducible component is defined by a mixed vanishing of prime factors of $z_j$ and $w^i_j$. We conclude that
every irreducible component admits at worst normal crossing
singularities and they meet in normal crossing way.  A priori, we cannot exclude the case when some prime divisor
of irreducible component has self-intersection. If in addition, we know that every prime divisor of irreducible component is smooth,
then we can conclude that every irreducible component is smooth. This is the case when $g=1$, but not the case when $g=2$, hence
it should not be the case for all $g \ge 2$.

When $g=1$, the matrix $\Phi=(c_{ij} \prod_{q \in N_{[ a_i,\delta_j]} } \zeta_q)$ in the local equation in \eqref{localEq} has only one row.
We can absorb all the invertible coefficients $c_{ij}$ into the frames. Then, a simple analysis implies that all the singularities
are caused by the presence of nodes. We call these singularities of topological type.
These are straightforward to resolve as in Theorem \ref{thm:m1}.

When $g=2$,  the matrix $\Phi=(c_{ij} \prod_{q \in N_{[ a_i,\delta_j]} } \zeta_q)$ in the local equation in \eqref{localEq} has two rows.
We classify the singularities into two types. Topological type: the ones caused solely by the presence of nodes. Geometric type:
 the ones caused by vanishing of minor determinants of the matrix $(c_{ij})$.  This classification works for all genera.
 In genus two, we prove that the singularities of
 geometric type means: Weierstrass and conjugate rational tails,
 double cover of rational curves (hyperelliptic). Thus, one may also call the geometric type
 the Brill-Noether type. These are not too hard to handle when $g=2$.  In Theorems \ref{thm:m2}
 and \ref{thm:m3}, the first two rounds of sequences of smooth blowups are along centers of topological type. Roughly, the first round deals with
 $(1\times 1)$ minors; the second round deals with $(2\times 2)$ minors. Both rounds  ignore singularities of geometric type. Then the third rounds
 resolve singularities of Weierstrass and conjugate rational tails. The fourth round resolves singularities caused by hyperelliptic maps.

 For general $g$,  it is not hard for the author
  to $``$see$"$ how to resolve singularities of topological type. The remainder singularities of geometric type
 are fortunately of high codimension. I believe that a canonical partial relative resolution that resolves singularities of topological type
 will be sufficient to prove Conjecture \ref{conj} (the details will be in a forthcoming publication).

%\section{Partial relative resolution}
  %We should mention that although $\fM_g^{\rm wt}$ is sufficient to treat all singularities cause by nodes
 %(thus, of topological type) and is easy to use, it is not sufficient to deal with all other singularities (of geometry type,
 %i.e., anything of Brill-Noether flavor) as soon as $g>1$.

\section{Applications to the GW theory }\lab{applications}
In this section, we simply  summarize the approach of \cite{CL12}.

The GSW invariants originate from Guffin and Sharpe's  work  [GuSh08].  Let
$\overline{M}_1(\PP^4,d)^p$ be the DM stack of pairs
$([u,C], p)$ with  $[u,C] \in \overline{M}_1(\PP^4,d)$ and  $p \in H^0(C, u^*\sO_{\Pf}(-5)\otimes \omega_C)$.
This is the moduli of stale maps with $p$-fields.  Chang and Li [CL11] showed that it admits with a perfect obstruction theory.
We shall  use the following shorthands
$$\cX=\overline{M}_1(\PP^4,d), \quad \cY=\overline{M}_1(\PP^4,d)^p.$$ The
polynomial $x_1^5+\ldots+x_5^5$ induces a cosection of its obstruction sheaf
%\beq\label{si0}
$\sigma : {\cal O}b_{\cY}\lra \sO_{\cY}$
%\eeq
whose  degeneracy locus (where $\sigma$ fails to be surjective) is
$$\overline{M}_1(Q,d)\sub \cY, \quad Q=(x_1^5+\ldots+x_5^5=0)\sub \Pf,
$$
which is proper.  The cosection localized virtual class  of Kiem-Li ([KL10])
defines a  localized virtual cycle
$[\cY]\virt_\sigma\in A_0 \overline{M}_1(Q,d).
$
The GSW-invariant of the quintic $Q$ is
\beq \lab{p-Q} N_{1,d}^p = \deg [\cY]\virt_\sigma.
\eeq

\begin{theo}  {\rm ([CL11])}
The GSW-invariant coincides with the GW-invariant of the quintic $Q$ up to a sign:
$$N_{1,d}^p=(-1)^{5d}\cdot N_{1,d}
$$ where $N_{1,d}$ is the (usual) genus $1$, degree $d$ GW invariant of the smooth quintic $Q \sub \Pf$.
\end{theo}

This theorem translates the problem from over $\cX$ to over $\cY$. To prove Conjecture \ref{conj},
it requires to study the separation of the virtual cycle $[\cY]\virt_\sigma$.  For this, we need some structural results
on its intrinsic normal cone. So,
we let $(f_\cX, \pi_\cX): \cC_\cX  \lra \Pf \times \cX$ be the universal family of $\cX$; let
$$(f_\cY, \pi_\cY): \cC_\cY  \lra \Pf \times \cY, \;\; \psi_\cY \in \Gamma (\cY, \cP_\cY), \;\; \cP_\cY = f_\cY^* \sO(-5) \otimes \omega_{\cC_\cY/\cY}$$
be the universal family of $\cY$.
A  perfect obstruction theory of $\cX$ relative to  the Artin stack $\fM$ of  genus $g$ nodal curves  is ([BF97])
$$\phi_{\cX/\fM}: ( E_{\cX/\fM})\dual  \lra L^\bullet_{\cX/\fM},\quad  E_{\cX/\fM}\defeq
R^\bullet\pi_{\cX\ast} f_{\cX}^*T_{\Pf}. $$
Using the Euler sequence,
  a  perfect obstruction theory of $\cX$ relative to $\fP_1$ is (cf. [CL11])
$$ \phi_{\cX/\fP_1}: ( E_{\cX/\fP_1})\dual  \lra L^\bullet_{\cX/\fP_1},\;  E_{\cX/\fP_1}\defeq
R^\bullet\pi_{\cX\ast} \cL_{\cX}^{\oplus 5},  \; \cL_\cX=f_\cX\sta\sO(1).$$
Chang   and Li also worked out
a perfect obstruction theory of $\cY$ relative to $\fP_1$:
$$
\phi_{\cY/\fP_1}:(E_{\cY/\fP_1})\dual
\lra L_{\cY/\fP_1}^\bullet  ,\; E_{\cY/\fP_1}\defeq   R^\bullet \pi_{\cY\ast}(\cL_\cY^{\oplus 5}\oplus \cP_\cY), \; \cL_\cY= f_\cY^* \sO(1).
$$
The cohomology sheaf
$$\Ob_{\cY/\fP_1}\defeq  H^1(E_{\cY/\fP_1})=R^1 \pi_{\cY\ast}(\sL_\cY^{\oplus 5}\oplus \cP_\cY) $$
is the relative obstruction sheaf of $\phi_{\cY/\fP_1}$.

The following theorem  will be useful.

\begin{theo} \lab{thm:p}
We let $\widetilde\fM_1^{\rm wt} \lra \fM_1^{\rm wt}$ be the successive blowup of $\fM_1^{\rm wt}$
 along  the smooth closed substacks $\Theta_2, \Theta_3, \cdots,$ and so on.
If we let $$\widetilde{M}_1 (\Pf,d)^p = \widetilde\fM_1^{\rm wt} \times_{ \fM_1^{\rm wt}}
\overline{M}_1(\Pf,d)^p, $$ then we have
\begin{enumerate}
%\item  the main component of $\widetilde{M}_1 (\Pn,d)$ is smooth;
\item  every (including the main) irreducible component of $\widetilde{M}_1 (\Pf,d)^p$ is smooth;
\item  the entire DM stack $\widetilde{M}_1 (\Pf,d)^p$ has at worse normal crossing singularities; further
\item the  object $R^\bullet \pi_{\cY\ast}(\cL_\cY^{\oplus 5}\oplus \cP_\cY)$, upon pulling back to $\widetilde{M}_1 (\Pf,d)^p $,
becomes locally diagonalizable.
\end{enumerate}
\end{theo}

This is analogous to Theorem \ref{thm:m1}  and can be proved parallely.

We let $\tcy$ be $\widetilde{M}_1 (\Pn,d)^p $ and $\tilde\fP_1=  \widetilde\fM_1^{\rm wt} \times_{ \fM_1^{\rm wt}}
\fP_1.$ Then we have a canonical morphism $\tcy  \to \tilde\fP_1$.
By working out the relative perfect obstruction theory of $\tcy\to\tilde\fP_1$, we obtain its obstruction complex $E_{\tcy/\tilde\fP_1}$.
The intrinsic normal cone $\bC_{\cpx/\tilde{\fP}_1}$ of $\cpx\to \tilde{\fP}_1$ is embedded in $h^1/h^0(E_{\cpx/\tilde{\fP}_1})$.

 The blowup $\tilde\cX$ is a union of smooth DM stacks: one, denoted  $\tilde\cX'$,
 is the proper transform of the main component $\cX' $ of $\cX$, the rest are indexed by partitions $\mu$   of $d$.
 In formula, $$\tilde\cX = \tilde\cX' \cup (\cup_{\mu \vdash d} \tilde\cX_\mu).$$
 Geometrically, general points of  $\tilde\cX'$ are stable morphisms with smooth domain curves;
general points of $\tilde\cX_\mu$ lie over stable morphisms whose domain curves are the union of a smooth genus one curve and
 connected trees of rational curves such that the genus one curve is contracted and the degrees of the stable morphism on the connected trees
 form the partition $\mu$. The blowup $\tcy$ has the similar decomposition
 $$\tilde\cY = \tilde\cY'  \cup (\cup_{\mu \vdash d} \tilde\cY_\mu).$$

 Theorem \ref{thm:p} implies that

\begin{lemm}\lab{normalConesAre}
 Let $\bC_{\tcy/\tilde\fP_1}\sub h^1/h^0(E_{\tcy/\tilde\fP_1})$ be the intrinsic normal cone embedded in
 $h^1/h^0(E_{\cpx/\tilde{\fP}_1})$ via the obstruction
theory $\phi_{\tcy/\tilde\fP_1}$.  Then
\begin{enumerate}
\item away from $\cup_{\mu \vdash d} \tilde\cY_\mu$, $\bC_{\tcy/\tilde\fP_1}$ is the zero section of   $h^1/h^0(E_{\tcy/\tilde\fP_1})$;
\item away from $ \tilde\cY' $,  it is a rank two subbundle stack of  $h^1/h^0(E_{\tcy/\tilde\fP_1})$.
\end{enumerate}
\end{lemm}

Lemma \ref{normalConesAre} allows us to conclude that
the cone $\bC_{\cpx/\tfd}$ can be separated as
$$[\bC_{\cpx/\tfd}] =[\bC']+\sum_{\mu\vdash d}\,  [\bC_\mu],
$$
where $\bC'$ is an irreducible cycle lying over $\tcy_{\rm pri}$;
each $\bC_\mu$ lies over $\cpx_\mu$.  (The cycles $\bC_\mu$ need not to be irreducible.)
Thus, applying Kiem-Li's cosection localized virtual class, we obtain
\beq\label{sep-vir}
\deg 0^!_{\tilde\sigma,\text{loc}}[\bC_{\cpx/\tfd}]  =0^!_{\tilde\sigma,\loc}[\bC']+
\sum_{\mu\vdash d} 0^!_{\tilde\sigma,\loc}[\bC_\mu],
\eeq
where  the cosection $\tilde\sigma:  {\cal O}b_{\tilde\cY}\lra \sO_{\tilde\cY}$ is induced from   $ \sigma : {\cal O}b_{\cY}\lra \sO_{\cY}$.

Now, a general standard argument implies that $$\deg 0_{\tilde\sigma,\loc}^![\bC']=(-1)^{5d} N_{1,d}'.$$
By some accurate dimensional arguments, one finds the vanishing
$$\deg 0_{\tilde\sigma,\loc}^![\bC_\mu]=0$$ for all  $\mu\ne (d)$ where
 $(d)$ be the partition of $d$ into a single part  (i.e., the non-partition of $d$).
Finally,  Chang and Li calculate to prove
$$\deg 0_{\tilde\sigma,\loc}^![\bC_{(d)}]=\frac{(-1)^{5d}}{12} N_{0,d}.$$
Putting these together, one gets $N_{1,d}=N_{1,d}'+ \frac{1}{12} N_{0,d}$.

%\section{Alternative Modular Compactifications}% of $M_g(\G\bR(r, \CC^n),d)$ }\lab{alternative}

\section{Modular Resolution Program}\lab{pmr}

Suppose we have  a moduli space $\fM$ and we  let $\widetilde\fM \to \fM$ be a blowup. In general, it is not clear how to provide  $\widetilde\fM$ a modular meaning.
There are a few interesting cases where we know that  the blowups have geometric meanings. Here we mention one class of such examples.

We let $\Gr(k,V)$ be the Grassmannian  of $k$-dimensional subspaces
in a vector space $V \cong \kk^n$ where $\kk$ is the base field  $(0<k<n$).
Consider  the variety $\cQ_d$ of degree $d$ maps from ${\bf P}^1$ to  $\Gr(k,V$).
It comes with a naive compactification, the Grothendieck Quot schome $\overline{\cQ}_d$
(= Quot$^{d,n-k}_{V_{{\bf P}^1/{\bf P}^1/{\bf k}}}$).
The compactification is smooth but the boundary $\overline{\cQ}_d \setminus \cQ_d$ has rather intricate singularities.
 It comes equipped with a natural filtration  by  closed subsets
 $$Z_{d,0} \subset Z_{d,1} \subset \cdots \subset Z_{d,d-1}=\overline{\cQ}_d \setminus \cQ_d$$
 where $Z_{d,r}$ consists of non-locally free sheaves whose torsion parts have degree at least $d-r$.

 \begin{theo}\lab{hls} {\rm ([HLS11], [Shao11]))}   The singularities of $Z_{d,r}$ can be resolved by repeatedly blowing up
$Z_{d,0}$, $Z_{d,1}, \cdots, Z_{d,r-1}$. Consequently, by successively blowing-up the Quot scheme $\overline{\cQ}_d$
along  $Z_{d,0}$, $Z_{d,1}, \cdots, Z_{d,d-1}$, we obtain a smooth compactification
$\widetilde{\cQ}_d$   such that the boundary $\tilde{\cQ}_d \setminus
 \cQ_d$ is a  simple normal crossing divisor.
 \end{theo}

Next, it is natural to ask whether
$\widetilde{\cal Q}_d$  admits any modular interpretation. We provide an
affirmative answer to this question.  For any coherent sheaf $\cE$ over $\Po$,
we let $\cE_t$ denote the torsion subsheaf of $\cE$
and $\cE_f = \cE/\cE_t$ denote the locally free part of $\cE$.

\begin{defi} \lab{def:cq} A complete quotient of $V_{\Po}$ consists of
$(\cE^0, \cdots, \cE^\ell)$ where $\cE^0$ is a coherent quotient sheaf of $V_{\Po}$,  and for every $1 \le i \le \ell$,
$\cE^i$ is a {\it nonsplit} extension of $\cE^{i-1}_t$ by $\cE^{i-1}_f$
such that the last sheaf $\cE^\ell$  is the unique one that is locally free. (We allow $\ell=0$.)
%The Hilbert polynomial of the complete quotient is defined to the Hilbert polynomial of its initial sheaf $\cE^0$.
 \end{defi}

\begin{theo}\lab{hs}
{\rm ([HS13])}  %Fix the Grassmannian $\Gr(k,V)$,  $1 \le k \le n-1$.
The variety $\widetilde{\cal Q}_d$in Theorem \ref{hls} parameterizes {\it complete quotients} of $V_{\Po}$
%$\sO_{\Po}^{\oplus n}$
of degree $d$ and rank $n-k$.
\end{theo}

%The blowup morphism   $\widetilde{\cal Q}_d \to \overline{\cal Q}_d$ is simply the forgetting map.
There are further interesting  questions.
%(1) generalize to  curves of higher genus; (2) apply to enumerative problems.
For $g=1$, the analogous of $\overline{\cal Q}_d$  is MOP's moduli $\overline{\rm SQ}_1(\Gr(r, V),d)$ of stable quotients.
A smooth blowup  $\widetilde{\rm SQ}_1(\Gr(r, V),d) \to \overline{\rm SQ}_1(\Gr(r, V),d)$,
{\it resolving all the boundary singularities},  is currently being worked out by
Thomas  D.  Maienschein %in his Ph.D thesis at the University of Arizona
(\cite{tdm}).  The above should generalize to all  high genera.

\medskip
The author believes the following.

\begin{conj} For every $g>0$, there is a smooth moduli stack $\fM$ admitting a dominating morphism $\bMPdg \to \fM$;
further, there is also another smooth {\bf moduli} stack $\widetilde\fM$ admitting a birational dominating morphism to  $\fM$ such that if we let
$$
\begin{CD}
\widetilde{M}_g(\Pn,d)= \widetilde\fM \times_{\fM} \overline{M}_g(\Pn,d),  \end{CD}
 $$
then  $\widetilde{M}_g(\Pn,d)$ inherits moduli meaning. Further, suppose $d > 2g-2$, then we have
\begin{enumerate}
\item  the main component of $\widetilde{M}_g(\Pn,d)$ is smooth;
\item  the entire stack $\widetilde{M}_g(\Pn,d)$ has at worst normal crossing singularities.
\end{enumerate}
\end{conj}

This is a kind of reformulation of Conjecture \ref{r-r}.   Here, the smoothness can subject to re-definition. For $g=1$, it is believed
that $\widetilde\fM$ can be taken to be the moduli of weighted  genus one curves with certain admissible log structures
(this is due to Qile Chen).

\medskip

%\medskip
Beyond  reach, the author ends this note  with an unlimitedly wild speculation.

\begin{conj}
Let $\fX$ be a singular moduli stack with the main component $\fX'$.
Then under some reasonable hypotheses on $\fX$,
there is  another moduli stack $\widetilde\fX$
with a smooth main component $\widetilde\fX'$,  at worst normal crossing singularities, and
admitting a dominant morphism to $\fX $ such that it restricts  to a birational morphism $\widetilde\fX' \to \fX'$. \end{conj}

Derived from {\it Anna Karenina},
$$\hbox{\it Regularities are all alike; every singularity is singular in its own way.}$$
The conjecture  would say that every singularity can be resolved because there are geometric reasons to let us achieve it (not because we can implant an algorithm). de Jong's alteration has this flavor.
We hope to furnish more interesting examples of this in the future publications.

\bigskip

%\centerline{\bf REFERENCES}

\end{document}